\documentclass{article}
\usepackage[hyphens]{url}
\usepackage{amsmath,amsthm,color,bm}

\newenvironment{newreferences}
               {\section*{References}
                \raggedright
                \begin{list}{}{\setlength{\itemsep}{0pt}
                               \setlength{\parsep}{0pt}
                               \setlength{\labelwidth}{0pt}
                               \setlength{\leftmargin}{12pt}
                               \setlength{\labelsep}{0pt}}
                \setlength{\itemindent}{-12pt}
               }{\end{list}}

\def\sp{$\spadesuit$}
\def\cl{$\clubsuit$}     
\def\di{$\diamondsuit$}
\def\he{$\heartsuit$}

\def\J{{\rm J}}        
\def\Q{{\rm Q}}        
\def\K{{\rm K}}        
\def\A{{\rm A}}        

\def\P{{\rm P}}        
\def\E{{\rm E}}        
\def\Var{{\rm Var}}
\def\Cov{{\rm Cov}}

\theoremstyle{plain}

\theoremstyle{remark}

\theoremstyle{definition}
\newtheorem{example}{Example}

\allowdisplaybreaks

\begin{document}
\title{Variance reduction \\ in Texas hold'em and in video poker}
\author{Stewart N. Ethier\thanks{Department of Mathematics, University of Utah. Email: ethier@math.utah.edu.}}
\date{}
\maketitle

\begin{abstract}
In Texas hold'em, after an all-in bet is made and called before the flop, the turn, or the river, the two players sometimes agree to run it $n$ times, meaning that the remaining five, two, or one cards are dealt out not just once but $n$ times successively without replacement, with $1/n$ of the pot attached to each run.  In $n$-play video poker, five cards are dealt exactly as in the conventional single-play game.  After the player chooses which cards to hold, new cards are drawn to replace the discards, not just once but $n$ times independently, with $1/n$ of the bet attached to each draw.  In both scenarios the players are attempting to reduce the variance of the return without changing the mean.  We quantify the extent to which the variance is reduced.
\end{abstract}

\section{Introduction}\label{intro}

The most important statistic of a wager is the mean return.  The second most important is arguably the variance of the return.  Sometimes a gambler has a choice between several wagers with the same mean return but different variances of return.  Should he try to minimize or maximize the variance?  If the gambler is a professional who has a long-term advantage, he would want to minimize the variance to make actual return closer to what is expected.  

We consider two scenarios where this occurs, one in Texas hold'em and the other in video poker.  In Texas hold'em, after an all-in bet is made and called before the flop, the turn, or the river, the two players sometimes agree to run it $n$ times, meaning that the remaining five, two, or one cards are dealt out not just once but $n$ times successively without replacement, with $1/n$ of the pot attached to each run.  Mean return is unaffected.  By how much is the variance reduced?

Running it $n$ times may appear to be a recent development in poker, but we have reliable evidence that it was used back in the 1960s in seven-card stud high-low games in New York City, often with $n=4$.  Malmuth (2021, p.~107) argues that running it twice should not be allowed by cardrooms because its effect is to make it more likely that the stronger players will win in the short term.  Fewer short-term wins for the weaker players might keep them from returning.  There is a delicate balance between luck and skill in poker, and it must be maintained for the health of the enterprise.

Another scenario is in $n$-play video poker.  Five cards are dealt exactly as in the single-play game.  After the player chooses which cards to hold, new cards are drawn to replace the discards, not just once but $n$ times independently, with $1/n$ of the bet attached to each draw.  Does the $n$-play game reduce the variance, relative to the single-play game?  The answer is yes, but it might not be as obvious as in the previous scenario.  In the well-regarded textbook, \textit{Practical Casino Math}, Hannum and Cabot (2005, p.~162) write,
\begin{quote}
Without delving further into the mathematics, suffice it to say that for multi-play video poker, there is no change in the expected value (house edge) but the variance (volatility) increases.
\end{quote}
For support they cite Jazbo (1999) and Shackleford (2012), who did delve further into the mathematics.  We will explain what we believe led to this mistaken conclusion from each of these authors.

The two scenarios differ in two important respects, so must be treated separately.  In the Texas hold'em setting, agreement to run it $n$ times is reached only near the end, after an all-in bet is made and called, whereas, in the video poker setting, agreement to play it $n$ times is reached at the beginning.\footnote{This suggests an interesting possibility: Video poker machines could be redesigned to allow the player to decide, based on his initial hand, whether to play it $n$ times.  Currently, the player must make this decision before seeing his initial hand.} In the Texas hold'em setting the $n$ runs are not independent, whereas in the video poker setting the $n$ plays are independent, at least conditionally.

Section~\ref{sec:poker} treats the Texas hold'em problem, and Section~\ref{sec:videopoker} treats the video poker problem.  We summarize our conclusions in Section~\ref{conclusions}.

\section{Texas hold'em, run $\bm n$ times}\label{sec:poker}

We assume that the reader is familiar with the terminology of Texas hold'em poker.  It sometimes happens, in Texas hold'em that, after an all-in bet is made and called before the flop, the turn, or the river, one player suggests, and his opponent agrees, to run it $n$ times, meaning that the remaining five, two, or one cards are dealt out not just once but $n$ times successively, without replacement, from the original deck, with 1/n of the pot attached to each run.  For a specified player, let $R_1,R_2,\ldots,R_n$ be an exchangeable sequence with $R_i$ equaling 1 if the player wins, 0 if the player loses, and 1/2 if the player ties on the $i$th run.  Then 
\begin{equation}\label{samplemean-2}
\frac{R_1+R_2+\cdots+R_n}{n}
\end{equation}
represents the portion of the pot that the specified player receives at the conclusion of the hand.  We do not consider three or more players in an all-in situation because that involves the main pot and a side pot, and we want to avoid such complications.  

The mean and variance of \eqref{samplemean-2} are, by exchangeability,
\begin{equation}\label{mean-formula}
\E\bigg[\frac{R_1+R_2+\cdots+R_n}{n}\bigg]=\frac{1}{n}\sum_{i=1}^n\E[R_i]=\E[R_1]
\end{equation}
and
\begin{align}\label{variance-formula}
\Var\bigg(\frac{R_1+R_2+\cdots+R_n}{n}\bigg)&=\frac{1}{n^2}\sum_{i=1}^n\Var(R_i)+\frac{1}{n^2}\,2\sum_{1\le i<j\le n}\text{Cov}(R_i,R_j)\nonumber\\
&=\frac{1}{n}\,\Var(R_1)+\bigg(1-\frac{1}{n}\bigg)\Cov(R_1,R_2).
\end{align}
Thus, the mean return for running it $n$ times is the same as for running it once, while the variance for running it $n$ times is less than $1/n$ times the variance for running it once.  This relies on a heuristically clear but unproved observation, namely that the covariance term is negative, the result of sampling without replacement.  Heuristically, a win on the first run reduces the chance of a win on the second run because the cards used to achieve the win on the first run are unavailable for the second run.

\subsection{One card to come}

The simplest case to analyze is when the all-in bet and call occur before the river.  Let the specified player be the one with the weaker hand.  Typically, among the 44 unseen cards, there are $o$ ``outs,'' that is, $o$ cards that will give the underdog the win; the remaining $44-o$ cards will preserve the win for the favorite.  Notice that $R_1+R_2+\cdots+R_n$ has a hypergeometric distribution, hence
\begin{equation}\label{hypergeo-1}
\E\bigg[\frac{R_1+R_2+\cdots+R_n}{n}\bigg]=\frac{o}{44}
\end{equation}
and
\begin{equation}\label{hypergeo-2}
\Var\bigg(\frac{R_1+R_2+\cdots+R_n}{n}\bigg)=\frac{1}{n}\,\frac{o}{44}\bigg(1-\frac{o}{44}\bigg)\bigg(1-\frac{n-1}{43}\bigg),
\end{equation}
results that are consistent with \eqref{mean-formula} and \eqref{variance-formula}.  Here $1\le n\le 44$ is implicitly assumed, and with $n=44$ the variance would be 0.  But in practice $n$ is further restricted by the number $m$ of hands that have been mucked and by the convention of burning a card before the flop, before the turn, and before the river.  Thus $1\le n\le 44-(2m+3)$.  The last factor on the right side of \eqref{hypergeo-2} is called the \textit{finite-population correction factor}.

\begin{example}\label{ex1a}
Player 1 has 10\sp-10\cl, Player 2 has \K\he-\K\di, the flop is \K\sp-\Q\sp-\J\sp, and the turn card is 7\di.  An all-in bet is made and called before the river.  Player 2 is the favorite with three kings, while Player 1 has an open-ended straight-flush draw.  Player 1 has 15 outs, and there is no possibility of a tie.  See Table~\ref{tab:1}.

\begin{table}[htb]
\caption{In Example~\ref{ex1a}, the mean, variance, and standard deviation of Player 1's return, in units of pot size, when the players agree to run it $n$ times.  From \eqref{hypergeo-1} and \eqref{hypergeo-2} with $o=15$.\label{tab:1}}
\catcode`@=\active \def@{\phantom{0}}
\begin{center}
\begin{tabular}{rccc}
    & &          & standard \\
$n$ & @@mean@@ & variance & deviation \\
\noalign{\smallskip}\hline\noalign{\smallskip}
1 & 0.340909 & 0.224690 & 0.474015 \\
2 & 0.340909 & 0.109732 & 0.331259 \\
3 & 0.340909 & 0.071413 & 0.267232 \\
4 & 0.340909 & 0.052254 & 0.228590 \\
\noalign{\smallskip}\hline
\end{tabular}
\end{center}
\end{table}

\end{example}

We can generalize this result slightly, allowing for ties.  In addition to the $o$ outs for the underdog, suppose there are $t$ cards that will result in a tie; the remaining $44-o-t$ cards will preserve the win for the favorite.  We can use \eqref{mean-formula} and \eqref{variance-formula} to find that
\begin{equation}\label{mean-outs/ties}
\E\bigg[\frac{R_1+R_2+\cdots+R_n}{n}\bigg]=\frac{o+\frac12 t}{44}
\end{equation}
and
\begin{align}\label{var-outs/ties}
&\Var\bigg(\frac{R_1+R_2+\cdots+R_n}{n}\bigg)\nonumber\\
&\quad{}=\frac{1}{n}\bigg(\frac{o+\frac14 t}{44}-\bigg(\frac{o+\frac12 t}{44}\bigg)^2\bigg)\nonumber\\
&\qquad{}+\bigg(1-\frac{1}{n}\bigg)\bigg(\frac{o(o-1)+\frac12 2ot+\frac14 t(t-1)}{44(43)}-\bigg(\frac{o+\frac12 t}{44}\bigg)^2\bigg).
\end{align}
The latter formula does not have an elegant simplification as in \eqref{hypergeo-2}.

\begin{example}\label{ex1b}
Player 1 has 10\di-10\cl, Player 2 has \K\he-\K\di, the flop is \K\sp-\Q\sp-\J\sp, and the turn card is A\sp.  An all-in bet is made and called before the river.  Player 1 is the favorite with an ace-high straight, while Player 2 has three kings.  Player 2 can win with a jack, queen, king, or ace, and can tie with a spade.  Thus, $o=10$ and $t=9$ in \eqref{mean-outs/ties} and \eqref{var-outs/ties}, which gives 
\begin{equation}\label{mean-outs10,ties9}
\E\bigg[\frac{R_1+R_2+\cdots+R_n}{n}\bigg]=\frac{29}{88}
\end{equation}
and
\begin{equation}\label{var-outs10,ties9}
\Var\bigg(\frac{R_1+R_2+\cdots+R_n}{n}\bigg)=\frac{1}{n}\,\frac{1315}{7744}-\bigg(1-\frac{1}{n}\bigg)\frac{1315}{332992}.
\end{equation}
See Table~\ref{tab:2}.

\begin{table}[htb]
\caption{In Example~\ref{ex1b}, the mean, variance, and standard deviation of Player 2's return, in units of pot size, when the players agree to run it $n$ times.  From \eqref{mean-outs10,ties9} and \eqref{var-outs10,ties9}.\label{tab:2}}
\catcode`@=\active \def@{\phantom{0}}
\begin{center}
\begin{tabular}{rccc}
    & &          & standard \\
$n$ & @@mean@@ & variance & deviation \\
\noalign{\smallskip}\hline\noalign{\smallskip}
1 & 0.329545 & 0.169809 & 0.412079 \\
2 & 0.329545 & 0.082930 & 0.287976 \\
3 & 0.329545 & 0.053970 & 0.232315 \\
4 & 0.329545 & 0.039490 & 0.198722 \\
\noalign{\smallskip}\hline
\end{tabular}
\end{center}
\end{table}

\end{example}

\subsection{Two cards to come}

A more complicated case to analyze is when the all-in bet and call occur before the turn.  When they run it twice, the number of ways for the hand to be played out is
\begin{equation*}
\binom{45}{2}\binom{43}{2}=(990)(903)=893970,
\end{equation*}
less than 900 thousand, so an app can be written, which will run almost instantly.  If $m$ is the number of hands that have been mucked, the constraint on the number $n$ of runs is $1\le n\le (45-(2m+3))/2=21-m$.  But to give some insight into what is involved, we consider an example.

\begin{example}\label{ex2}
This hand occurred in a televised game, circa 2013; a partial video can be found at \url{https://www.dailymotion.com/video/x10aw9l}.  Player 1 (Ernest Wiggins) has \K\sp-\K\cl.  Player 2 (Phil Hellmuth) has \A\he-9\di.  The flop is 9\he-10\sp-9\sp.  Player 2 makes an all-in bet and Player 1 calls.  The pot has \$199{,}600.  Player 2 is the favorite, he suggests running it four times, and Player 1 agrees.  Run 1 results in \J\he-\A\di, a win for Player 2.  Run 2 results in \K\he-3\cl, a win for Player 1.  Run 3 results in 7\sp-8\sp, a win for Player 1.  Run 4 results in \Q\cl-\K\di, a win for Player 1.  The underdog wins three out of four.

Let us show how to analyze this example, the principles of which should generalize.  There is no chance of a tie, so it suffices to find the probability that Player 1 wins.  This event can occur in any of six ways:  (a) With a straight flush by drawing \Q\sp-\J\sp. (b) With four of a kind by drawing \K\he-\K\di. (c) With kings full of nines by drawing a red king and any other card except a red king and a 9\cl.  (d) With tens full of kings by drawing two of the three available tens.  (e) With a flush by drawing two of nine spades (ace excluded, king, 10, 9 unavailable), except for \Q\sp-\J\sp. (f) With a straight by drawing a queen and a jack, except \Q\sp-\J\sp.

To derive a formula for the probability of a Player 1 win that applies at both the first run and the second, we need to partition the deck into subsets, including aces, kings, queens, jacks, 10s, 9s, and others, both in spades and in nonspades.  So there are 11 subsets (see Table~\ref{subsets}).  We let $\bm m=(m_1,m_2,\ldots,m_{11})$ be the numbers of cards in these subsets.  Initially, $\bm m=(1,2,2,1,3,1,3,3,1,7,21)$.  The sum of these counts is 45.

\begin{table}[htb]
\caption{Subsets of cards in Example~\ref{ex2}.\label{subsets}}
\begin{center}
\begin{tabular}{lll}
& spades & nonspades \\
\noalign{\smallskip}\hline\noalign{\smallskip}
aces    & 1. \A\sp & 2. \A\cl, \A\di \\
kings   &          & 3. \K\di, \K\he \\
queens  & 4. \Q\sp & 5. \Q\cl, \Q\di, \Q\he \\
jacks   & 6. \J\sp & 7. \J\cl, \J\di, \J\he \\
10s     &          & 8. 10\cl, 10\di, 10\he \\
9s      &          & 9. 9\cl \\
2s--8s  & 10. 2\sp--8\sp & 11. 2--8 nonspades \\
\noalign{\smallskip}\hline
\end{tabular}
\end{center}
\end{table}

The probability that Player 1 wins, as a function of $\bm m$, is
\begin{align*}
&p(\bm m)\\
&\quad:=\P(R_1=1)\\
&\quad\;=\bigg(m_4 m_6 + \binom{m_3}{2} + m_1 (m-m_3-m_9)+ \binom{m_8}{2}+ \binom{m_4 + m_6 + m_{10}}{2}\\
&\quad\qquad{}  - m_4 m_6 + (m_4 + m_5) (m_6 + m_7) - m_4 m_6\bigg)\bigg/\binom{m}{2},
\end{align*}
where $m:=m_1 + m_2 + m_3 + m_4 + m_5 + m_6 + m_7 + m_8 + m_9 + m_{10} + m_{11}$.
We find that $p(1,2,2,1,3,1,3,3,1,7,21)=139/990$.\footnote{The careful reader will notice a discrepancy between our figure for Player 1's mean return (14\,\%) and the one in the video cited above (15\,\%).  Our figure is from the perspective of Players 1 and 2, who have seen only seven cards when they make their agreement to run it four times.  The video's figure is from the perspective of the television audience, who, thanks to hole-card cameras, have seen 15 cards at that time, including the eight cards mucked by the four players who are no longer active in the hand.  The eight additional cards are 5\cl-2\he, 6\he-5\he, 6\sp-3\sp, and J\sp-10\he, so the relevant probability is $p(1,2,2,1,3,0,3,2,1,5,17)=97/666$.

It should be acknowledged that our figure for Player 1's mean return (from the perspective of Players 1 and 2) relies on a questionable assumption, namely that all $\binom{m}{2}$ possible sets of turn and river cards are equally likely to appear.  We know that four hands were mucked, which provides some, albeit very limited, information about them.  For example, it was a near certainty (before it was confirmed) that none of these four mucked hands was a pair of aces, kings, queens, or jacks.}

In two runs, the probability that Player 1 wins twice is
\begin{align*}
&\P(R_1=1,\,R_2=1)\\
&=\bigg(m_4 m_6\, p(\bm m-\bm e_4-\bm e_6) + \binom{m_3}{2} p(\bm m-2\bm e_3)+ m_3 m_1\, p(\bm m-\bm e_3-\bm e_1)\\ 
&\quad{} + m_3 m_2\, p(\bm m-\bm e_3-\bm e_2)+ m_3 m_4\, p(\bm m-\bm e_3-\bm e_4) + m_3 m_5\, p(\bm m-\bm e_3-\bm e_5)\\
&\quad{}+ m_3 m_6\, p(\bm m-\bm e_3-\bm e_6) + m_3 m_7\, p(\bm m-\bm e_3-\bm e_7)+ m_3 m_8\, p(\bm m-\bm e_3-\bm e_8)\\ 
&\quad{}+ m_3 m_{10}\,p(\bm m-\bm e_3-\bm e_{10})+ m_3 m_{11}\,p(\bm m-\bm e_3-\bm e_{11})+ \binom{m_{10}}{2} p(\bm m-2\bm e_{10})\\
&\quad{}+ m_4 m_{10}\, p(\bm m-\bm e_4-\bm e_{10}) + m_6 m_{10}\, p(\bm m-\bm e_6-\bm e_{10})\\
&\quad{}+ \binom{m_{10}}{2} p(\bm m-2\bm e_{10}) + m_4 m_7\, p(\bm m-\bm e_4-\bm e_7)+ m_5 m_6\, p(\bm m-\bm e_5-\bm e_6)\\ 
&\quad{}+ m_5 m_7\, p(\bm m-\bm e_5-\bm e_7)\bigg)\bigg/\binom{m}{2},
\end{align*}
where $\bm m=(1,2,2,1,3,1,3,3,1,7,21)$, $\bm e_1=(1,0,0,0,0,0,0,0,0,0,0)$ for example, and $m$ is as above.  This equals $38/2365$.  

We conclude from \eqref{mean-formula} and \eqref{variance-formula} that
\begin{equation}\label{mean-Hellmuth}
\E\bigg[\frac{R_1+R_2+\cdots+R_n}{n}\bigg]=\frac{139}{990}
\end{equation}
and
\begin{equation}\label{Var-Hellmuth}
\Var\bigg(\frac{R_1+R_2+\cdots+R_n}{n}\bigg)=\frac{1}{n}\,\frac{118289}{980100}-\bigg(1-\frac{1}{n}\bigg)\frac{153643}{42144300}.
\end{equation}
See Table~\ref{tab:3}.

Recall that Player 1 won 3/4 of the pot.  That is about 3.68 standard deviations above the mean, a rather unusual event.

\begin{table}[htb]
\caption{In Example~\ref{ex2}, the mean, variance, and standard deviation of Player 1's return, in units of pot size, when the players agree to run it $n$ times.  From \eqref{mean-Hellmuth} and \eqref{Var-Hellmuth}.\label{tab:3}}
\catcode`@=\active \def@{\phantom{0}}
\begin{center}
\begin{tabular}{rccc}
    & &         & standard \\
$n$ & @@mean@@ & variance & deviation \\
\noalign{\smallskip}\hline\noalign{\smallskip}
1 & 0.140404 & 0.120691 & 0.347406 \\
2 & 0.140404 & 0.058523 & 0.241914 \\
3 & 0.140404 & 0.037800 & 0.194422 \\
4 & 0.140404 & 0.027438 & 0.165646 \\
\noalign{\smallskip}\hline
\end{tabular}
\end{center}
\end{table}

\end{example}

\subsection{Five cards to come}

A much more complicated case to analyze is when the all-in bet and call occur before the flop.  When they run it twice, the number of ways for the hand to be played out is
\begin{equation*}
\binom{48}{5}\binom{43}{5}=(1712304)(962598)=1648260405792,
\end{equation*}
more than 1.6 trillion, so an app can be written, but it may be too slow to be of practical value.  If $m$ is the number of hands that have been mucked, the constraint on the number $n$ of runs is $1\le n\le \lfloor(48-(2m+3))/5\rfloor$.  To give some insight into what is involved, we consider one of the few examples for which exact calculations are straightforward.

\begin{example}\label{ex5}
Player 1 has \A\di-\A\cl\ and Player 2 has \A\sp-\A\he.  Typically, the result is a tie and the pot is split (``chopped'').  The only way to have a different result is if the board shows four or more cards of the same suit, which gives the player with the ace of that suit a flush (or straight flush) and the win.  

Let $D$ denote the set of 4-tuples of nonnegative integers summing to 5, that is, the set of all $(i,j,k,l)$ with $i,j,k,l\ge0$ and $i+j+k+l=5$.  The element $(i,j,k,l)$ represents the numbers of clubs, diamonds, hearts, and spades on the 5-card board.  Notice that $D$ has $\binom{8}{3}=56$ elements.  For a single run, the distribution of suits on the board is
\begin{equation*}
P(i,j,k,l):=\frac{\binom{12}{i}\binom{12}{j}\binom{12}{k}\binom{12}{l}}{\binom{48}{5}},\qquad i,j,k,l\ge0,\; i+j+k+l=5,
\end{equation*}
so
\begin{align}\label{Rdist}
\P(R=1)&=\sum_{(i,j,k,l)\in D: \max(i,j)\ge4}P(i,j,k,l)=\frac{1}{46},\nonumber\\
\P(R=\textstyle{\frac12})&=\sum_{(i,j,k,l)\in D:\max(i,j,k,l)\le3}P(i,j,k,l)=\frac{44}{46},\\
\P(R=0)&=\sum_{(i,j,k,l)\in D: \max(k,l)\ge4}P(i,j,k,l)=\frac{1}{46}.\nonumber
\end{align}

For two runs, it becomes
\begin{equation*}
P(i_1,j_1,k_1,l_1;i_2,j_2,k_2,l_2):=P(i_1,j_1,k_1,l_1)\frac{\binom{12-i_1}{i_2}\binom{12-j_1}{j_2}\binom{12-k_1}{k_2}\binom{12-l_1}{l_2}}{\binom{43}{5}}
\end{equation*}
for $i_1,j_1,k_1,l_1\ge0$, $i_1+j_1+k_1+l_1=5$, $i_2,j_2,k_2,l_2\ge0$, and $i_2+j_2+k_2+l_2=5$.  
This allows us to evaluate the joint distribution of $R_1$ and $R_2$ when we run it twice.  For example, the first of nine probabilities has the form
\begin{equation*}
\P(R_1=1,\,R_2=1)=\sum_{\substack{((i_1,j_2,k_1,l_1),(i_2,j_2,k_2,l_2))\in D\times D: \\ \max(i_1,j_1)\ge4\text{ and }\max(i_2,j_2)\ge4}}P(i_1,j_1,k_1,l_1;i_2,j_2,k_2,l_2).
\end{equation*}
Numerically, the resulting joint distribution is given by
\begin{equation}\label{jointdist}
\bordermatrix{& R_2=1 & R_2=\frac12 & R_2=0 \cr
R_1=1 & 132841 & 7180272 & 227238 \cr
R_1=\frac12 & 7180272 & 317414900 & 7180272 \cr
R_1=0 & 227238 & 7180272 & 132841 \cr}\frac{1}{346856146}.
\end{equation}
From \eqref{Rdist} we find that $\E[R_1]=1/2$ and $\Var(R_1)=1/92$, and from \eqref{jointdist} we get $\Cov(R_1,R_2)=-94397/693712292$.  

We conclude from \eqref{mean-formula} and \eqref{variance-formula} that
\begin{equation}\label{mean-AAvAA}
\E\bigg[\frac{R_1+R_2+\cdots+R_n}{n}\bigg]=\frac{1}{2}
\end{equation}
and
\begin{equation}\label{Var-AAvAA}
\Var\bigg(\frac{R_1+R_2+\cdots+R_n}{n}\bigg)=\frac{1}{n}\,\frac{1}{92}-\bigg(1-\frac{1}{n}\bigg)\frac{94397}{693712292}.
\end{equation}
See Table~\ref{tab:4}.

\begin{table}[htb]
\caption{In Example~\ref{ex5}, the mean, variance, and standard deviation of Player 1's return, in units of pot size, when the players agree to run it $n$ times.  From \eqref{mean-AAvAA} and \eqref{Var-AAvAA}.\label{tab:4}}
\catcode`@=\active \def@{\phantom{0}}
\begin{center}
\begin{tabular}{rccc}
    &  &        & standard \\
$n$ & mean & variance & deviation \\
\noalign{\smallskip}\hline\noalign{\smallskip}
1 & 1/2 & 0.010870 & 0.104257 \\
2 & 1/2 & 0.005367 & 0.073258 \\
3 & 1/2 & 0.003532 & 0.059435 \\
4 & 1/2 & 0.002615 & 0.051140 \\
\noalign{\smallskip}\hline
\end{tabular}
\end{center}
\end{table}

\end{example}

\section{$\bm n$-play video poker}\label{sec:videopoker}

In $n$-play video poker, the player bets one unit. He then receives five cards face up on the screen, with each of the $\binom{52}{5}$ possible hands equally likely. For each card, the player must then decide whether to hold that card or not. Thus, there are $2^5$ ways to play the hand. If he holds $k$ cards, then the following occurs $n$ times with the results conditionally independent: He is dealt $5-k$ new cards, with each of the $\binom{47}{5-k}$ possibilities equally likely. The player then receives his payout for the $n$ hands, which depends on the payout schedule and assumes $1/n$ of a unit is bet on each of the $n$ hands.  How much is the variance reduced, relative to the single-play game?  The answer to this question depends on the details of the game.  

Let the random variable $R_1$ denote the return, per unit bet, in the single-play game, and let the random vector $\bm X$ denote the hand initially dealt to the player, or more precisely the equivalence class it belongs to.  (In games without wild cards, we call two hands \textit{equivalent} if they have the same five denominations and if the corresponding denominations have the same suits after a permutation of $(\clubsuit,\diamondsuit,\heartsuit,\spadesuit)$.) Then the law of total variance tells us that
\begin{align}\label{totalvariance}
\Var(R_1)&=\Var(\E[R_1\mid\bm X])+\E[\Var(R_1\mid\bm X)]\nonumber\\
&=V_\text{deal}+V_\text{draw}.
\end{align}
We refer to the first expression on the right as the variance explained by the deal, and the second expression as the variance explained by the draw.  Other terms that are sometimes used are the \textit{between} and the \textit{within} variance, or the \textit{explained} and the \textit{unexplained} variance.  But in our problem there are two time-separated sources of randomness, the deal and the draw, and our terminology recognizes that.

Next, to introduce our notation, suppose temporarily that $n$ units are bet on an $n$-play video poker hand, one unit on each play, and let $R_1,R_2,\ldots,R_n$ be the returns from each such play.  Then 
\begin{equation}\label{samplemean}
\frac{R_1+R_2+\cdots+R_n}{n}
\end{equation}
represents the return \textit{per unit bet} on $n$-play video poker.  (By dividing by $n$, we are implicitly assuming that $1/n$ of a unit is bet on each play.) The evaluation of the mean and variance of the sum in \eqref{samplemean} is straightforward:
\begin{align*}
\E[R_1+R_2+\cdots+R_n]&=\E[\E[R_1+R_2+\cdots+R_n\mid\bm X]]\\
&=\E[n\,\E[R_1\mid\bm X]]\\
&=n\,\E[\E[R_1\mid\bm X]]\\
&=n\,\E[R_1],
\end{align*}
and
\begin{align*}
&\Var(R_1+R_2+\cdots+R_n)\\
&\quad{}=\Var(\E[R_1+R_2+\cdots+R_n\mid\bm X])+\E[\Var(R_1+R_2+\cdots+R_n\mid\bm X)]\\
&\quad{}=\Var(n\,\E[R_1\mid\bm X])+\E[n\,\Var(R_1\mid\bm X)]\\
&\quad{}=n^2\,V_\text{deal}+n\,V_\text{draw},
\end{align*}
where we are using the fact that $R_1,R_2,\ldots,R_n$ are conditionally i.i.d.\ given $\bm X$.  We conclude that
\begin{equation}\label{ERbar}
\E\bigg[\frac{R_1+R_2+\cdots+R_n}{n}\bigg]=\E[R_1]
\end{equation}
and
\begin{equation}\label{VarRbar-1}
\Var\bigg(\frac{R_1+R_2+\cdots+R_n}{n}\bigg)=V_\text{deal}+\frac{1}{n}V_\text{draw}.
\end{equation}
Alternative expressions for the variance include
\begin{equation}\label{VarRbar-2}
\Var\bigg(\frac{R_1+R_2+\cdots+R_n}{n}\bigg)=\Var(R_1)-\bigg(1-\frac{1}{n}\bigg)V_\text{draw}
\end{equation}
and
\begin{equation}\label{VarRbar-3}
\Var\bigg(\frac{R_1+R_2+\cdots+R_n}{n}\bigg)=\frac{1}{n}\,\Var(R_1)+\bigg(1-\frac{1}{n}\bigg)V_\text{deal},
\end{equation}
both of which follow from \eqref{VarRbar-1} and \eqref{totalvariance}.  

Equation \eqref{VarRbar-1} seems the most transparent, given that the variance of a sample mean is the variance of an observation divided by the sample size:  The variance explained by the deal is left unchanged, while the variance explained by the draw is divided by $n$.  Equation \eqref{VarRbar-2}, which first appeared in Ethier (2010, Problem 17.10), is useful because it tells us exactly by how much the $n$-play game reduces the variance, compared to the single-play game.  Equation~\eqref{VarRbar-3} is also interesting because it tells us exactly by how much the $n$-play game increases the variance, compared to $n$ independent single-play games, each for $1/n$ of a unit.

There is an alternative approach to the variance calculation using only the exchangeability of $R_1,R_2,\ldots,R_n$, namely (from \eqref{variance-formula})
\begin{equation}\label{VarRbar-4}
\Var\bigg(\frac{R_1+R_2+\cdots+R_n}{n}\bigg)=\frac{1}{n}\,\Var(R_1)+\bigg(1-\frac{1}{n}\bigg)\Cov(R_1,R_2).
\end{equation}

By \eqref{VarRbar-3} and \eqref{VarRbar-4},
\begin{align}\label{VarHannum}
\frac{\Var(R_1+R_2+\cdots+R_n)}{n}&=\Var(R_1)+(n-1)V_\text{deal}\nonumber\\
&=\Var(R_1)+(n-1)\Cov(R_1,R_2),
\end{align}
which is increasing in $n$.  This is ``variance-of-return per unit bet'' instead of ``variance of return-per-unit-bet.''  We believe confusion between these two expressions is what led to the mistaken conclusions of Hannum and Cabot (2005), Jazbo (1999), and Shackleford (2012).  We hasten to add that \eqref{VarHannum} is perfectly correct as an equation, but it should not be interpreted as the appropriate notion of variance for $n$-play video poker.

We note in passing the curious formula from \eqref{VarHannum}, namely
\begin{equation}\label{identity}
\Var(\E[R_1\mid\bm X])=\Cov(R_1,R_2),
\end{equation}
which applies when $R_1$ and $R_2$ are exchangeable and conditionally i.i.d.~given $\bm X$.  In particular, the covariance is nonnegative.

\subsection{Jacks or Better}

Let us examine the case of (full-pay) Jacks or Better, probably the most studied variant of video poker.  The pay table is $(800,50,25,9,6,4,3,2,1,0)$, with the positive entries representing the return, per unit bet, for a royal flush, straight flush, four of a kind, full house, flush, straight, three of a kind, two pairs, and a pair of jacks or better, respectively.  This assumes a maximum bet, typically five times the minimum bet.  There are 134459 equivalence classes of initial hands.  The optimal strategy has been known for nearly 30 years (Dancer, 1996) and it is essentially unique, meaning that the choice of optimal strategy does not affect the distribution of $R_1$.  It is also well known (Ethier, 2010) that, under optimal play,
\begin{equation*}
\E[R_1]=\frac{1653526326983}{1661102543100}\approx0.995439
\end{equation*}
and 
\begin{equation}\label{var-JoB}
\Var(R_1)=\frac{53846098447064372932173011}{2759261658693287357610000}\approx19.514676.
\end{equation}
Exact fractions are provided here because combinatorial problems deserve exact answers.

A few years ago, Ethier, Kim, and Lee (2019), using C++ and \textit{Excel}, found the distribution of $\E[R_1\mid\bm X]$ explicitly, and it was fully specified in the 35-page appendix of the ArXiv version of their paper.  It is a random variable with 1153 distinct values, and its variance can be computed to be
\begin{equation*}
V_\text{deal}=\frac{602864541441854523450139}{306584628743698595290000}\approx1.966389.
\end{equation*}
Here (and in \eqref{var-JoB}) we used \textit{Mathematica} because \textit{Excel} keeps only 15 digits of precision.
The same authors did not evaluate the distribution of $\Var(R_1\mid\bm X)$.  It could certainly be done by rerunning their C++ program with a few extra lines of code.  Nevertheless, by \eqref{var-JoB} and \eqref{totalvariance}, we can evaluate its mean, and that is all we need here:
\begin{equation*}
V_\text{draw}=\frac{55891953982463387918}{3185037467325338625}\approx17.548288.
\end{equation*}
It is not surprising that the variance explained by the draw is significantly greater than the variance explained by the deal.  A partial justification for this is that a drawn royal flush is about 16 times as likely as a dealt royal flush.

Table~\ref{JoB-data} evaluates the mean, variance, and standard deviation, of return per unit bet, from the $n$-play game for several choices of $n$.  For example, 100-play Jacks or Better reduces the standard deviation, relative to single-play, by a factor of 3.02.

\begin{table}[htb]
\caption{In the $n$-play versions of Jacks or Better video poker, the mean, variance, and standard deviation, of return per unit bet, assuming a maximum bet.  From \eqref{ERbar} and \eqref{VarRbar-1}, \eqref{VarRbar-2}, or  \eqref{VarRbar-3}.\label{JoB-data}}
\catcode`@=\active \def@{\phantom{0}}
\begin{center}
\begin{tabular}{rccc}
    &        &           & standard  \\
$n$ & mean   & variance  & deviation \\
\noalign{\smallskip}\hline\noalign{\smallskip}
@@1 & 0.995439 & 19.514676 & 4.417542 \\
@@3 & 0.995439 & @7.815818 & 2.795678 \\
@@5 & 0.995439 & @5.476046 & 2.340095 \\
@10 & 0.995439 & @3.721217 & 1.929046 \\
@25 & 0.995439 & @2.668320 & 1.633499 \\
@50 & 0.995439 & @2.317354 & 1.522286 \\
100 & 0.995439 & @2.141872 & 1.463513 \\
\noalign{\smallskip}\hline
\end{tabular}
\end{center}
\end{table}

\subsection{Deuces Wild}

Let us also examine the case of (full-pay) Deuces Wild, in which the four deuces are wild.  The pay table is $(800,200,25,15,9,5,3,2,2,1,0)$, with the positive entries representing the returns, per unit bet, for natural royal flush, four deuces, wild royal flush, five of a kind, straight flush, four of a kind, full house, flush, straight, and three of a kind.  This too assumes a maximum bet, typically five times the minimum bet.  There are 102359 equivalence classes of initial hands.  (Two initial hands are \textit{equivalent} if they have the same five denominations and if the corresponding non-deuce denominations have the same suits after a permutation of (\cl,\di,\he,\sp).)  The optimal strategy has long been known, but is not essentially unique (Ethier, 2010, p.~563), so we adopt the minimum-variance optimal strategy.  It is also well known (Ethier, 2010) that, under the stated strategy,
\begin{equation*}
\E[R_1]=\frac{32187682693}{31944279675}\approx1.007620
\end{equation*}
and 
\begin{equation}\label{var-DW}
\Var(R_1)=\frac{342713803167529293506213}{13265681051410035373125}\approx25.834618.
\end{equation}

Similar methods to those of Ethier, Kim, and Lee (2019) yield the distribution of $\E[R_1\mid\bm X]$.  It is a random variable with 8903 distinct values (!), and its variance can be computed to be
\begin{equation*}
V_\text{deal}=\frac{27769926596703007394837}{8843787367606690248750}\approx3.140049.
\end{equation*}
It follows from \eqref{var-DW} and \eqref{totalvariance} that
\begin{equation*}
V_\text{draw}=\frac{144566104812712980751}{6370074934650677250}\approx22.694569.
\end{equation*}

Table~\ref{DW-data} evaluates the mean, variance, and standard deviation, of return per unit bet, from the $n$-play game for several choices of $n$.  For example, 100-play Deuces Wild reduces the standard deviation, relative to single-play, by a factor of 2.77.

\begin{table}[htb]
\caption{In the $n$-play versions of Deuces Wild video poker, the mean, variance, and standard deviation, of return per unit bet, assuming a maximum bet.  From \eqref{ERbar} and \eqref{VarRbar-1}, \eqref{VarRbar-2}, or  \eqref{VarRbar-3}.\label{DW-data}}
\catcode`@=\active \def@{\phantom{0}}
\begin{center}
\begin{tabular}{rccc}
    &        &           & standard  \\
$n$ & mean   & variance  & deviation \\
\noalign{\smallskip}\hline\noalign{\smallskip}
@@1 & 1.007620 & 25.834618 & 5.082777 \\
@@3 & 1.007620 & 10.704905 & 3.271835 \\
@@5 & 1.007620 & @7.678963 & 2.771094 \\
@10 & 1.007620 & @5.409506 & 2.325834 \\
@25 & 1.007620 & @4.047832 & 2.011922 \\
@50 & 1.007620 & @3.593940 & 1.895769 \\
100 & 1.007620 & @3.366995 & 1.834937 \\
\noalign{\smallskip}\hline
\end{tabular}
\end{center}
\end{table}

Kim (2012) has also studied Double Bonus and Joker Wild video poker games, but a careful analysis of the nonuniqueness of the optimal strategy remains to be done.  For this reason, we can evaluate $V_\text{deal}$, which does not depend on which optimal strategy is used, but not $\Var(R_1)$, which does.  Therefore, we do not consider these games here.

\section{Conclusions}\label{conclusions}

In Texas hold'em, a useful upper bound on the variance of return (with return measured in units of pot size) after an all-in bet is called and the players agree to run it $n$ times is
\begin{equation*}\label{bound}
\Var\bigg(\frac{R_1+R_2+\cdots+R_n}{n}\bigg)\le\frac{1}{n}\,\Var(R_1).
\end{equation*}
With two or five cards to come, an app is probably needed to evaluate $\Var(R_1)$, but then one might as well evaluate not just the upper bound but
\begin{equation*}
\Var\bigg(\frac{R_1+R_2+\cdots+R_n}{n}\bigg)=\frac{1}{n}\,\Var(R_1)+\bigg(1-\frac{1}{n}\bigg)\Cov(R_1,R_2).
\end{equation*}
With only one card to come and no possibility of a tie, we have the simple formula 
\begin{equation*}
\Var\bigg(\frac{R_1+R_2+\cdots+R_n}{n}\bigg)=\frac{1}{n}\,\frac{o}{44}\bigg(1-\frac{o}{44}\bigg)\bigg(1-\frac{n-1}{43}\bigg)
\end{equation*}
in terms of the number $o$ of outs.

In video poker, the variance of return-per-unit-bet, $\Var(R_1)$, in a single-play game played optimally is often well known (e.g., in Jacks or Better, it is about 19.514676).  As discussed in Section~\ref{sec:videopoker}, $\Var(R_1)=V_\text{deal}+V_\text{draw}$, where $V_\text{deal}$ is the variance explained by the deal and $V_\text{draw}$ is the variance explained by the draw (e.g., in Jacks or Better, $V_\text{deal}\approx1.966389$ and $V_\text{draw}\approx17.548288$).  The variance of return-per-unit-bet in $n$-play video poker played optimally is then
\begin{equation*}
\Var\bigg(\frac{R_1+R_2+\cdots+R_n}{n}\bigg)=V_\text{deal}+\frac{1}{n}\,V_\text{draw}.
\end{equation*}
There are several other equivalent expressions.

\begin{newreferences}

\item Dancer, B. (1996) 9-6 Jacks or Better Video Poker: A Complete ``How to Beat the Casino'' Discussion. Published by the author.  (Available at UNLV Lied Library Special Collections.)

\item Ethier, S. N. (2010) \textit{The Doctrine of Chances: Probabilistic Aspects of Gambling}. Springer, Berlin and Heidelberg.  \url{http://www.math.utah.edu/~ethier/sample.pdf}.

\item Ethier, S. N., Kim, J. J., and Lee, J. (2019) Optimal conditional expectation at the video poker game Jacks or Better.   \textit{UNLV Gaming Research \& Review Journal} \textbf{23} 1--18.  \url{https://arxiv.org/abs/1602.04171}.

\item Hannum, R. C. and Cabot, A. N. (2005) \textit{Practical Casino Math}, Second Edition.  Institute for the Study of Gambling and Commercial Gaming, Reno, NV.  Trace Publications, Las Vegas, NV.

\item Jazbo (1999) An analysis of $N$-play video poker. \url{http://www.jazbo.com/videopoker/nplay.html}.

\item Kim, J. J. (2012) Optimal Strategy Hand-rank Table for Jacks or Better, Double Bonus, and Joker Wild Video Poker.
Masters thesis, McMaster University. \url{https://macsphere.mcmaster.ca/bitstream/11375/11863/1/fulltext.pdf}.

\item Malmuth, M. (2021) \textit{Cardrooms: Everything Bad and How to Make Them Better. An Analysis of Those Areas Where Poker Rooms Need Improvement}. Two Plus Two Publishing, Henderson, NV.

\item Shackleford, M. (2012) Standard deviation for multihand video poker. Wizard of Odds, \url{https://wizardofodds.com/games/video-poker/appendix/3/}.

\end{newreferences}

\end{document}